\documentclass[a4paper,10pt]{amsart}
\usepackage[german,english]{babel}
\usepackage[latin1]{inputenc}
\usepackage{ucs}
\usepackage{amsmath}
\usepackage{amsfonts}
\usepackage{amssymb}
\usepackage{amsthm}
\usepackage{fontenc}
\usepackage{graphicx}
\usepackage{url}
\usepackage{color}
\usepackage[T1]{fontenc}
\usepackage{a4wide}
\usepackage{booktabs}
\usepackage{float}
\usepackage{listings}

\newtheorem{remark}{Remark}

\newtheorem{conclusion}{Conclusion}

\newcommand{\bsx}{\boldsymbol{x}}
\newcommand{\bsy}{\boldsymbol{y}}
\newcommand{\bsz}{\boldsymbol{z}}
\newcommand{\bsn}{\boldsymbol{n}}
\newcommand{\bss}{\boldsymbol{s}}

\newcommand{\bsK}{\boldsymbol{K}}

\newcommand{\bsA}{\boldsymbol{A}}
\newcommand{\bsB}{\boldsymbol{B}}
\newcommand{\bsC}{\boldsymbol{C}}
\newcommand{\bsM}{\boldsymbol{M}}
\newcommand{\bsW}{\boldsymbol{W}}
\newcommand{\bsV}{\boldsymbol{V}}
\newcommand{\bsI}{\boldsymbol{I}}
\newcommand{\bsH}{\boldsymbol{H}}
\newcommand{\bsL}{\boldsymbol{L}}
\newcommand{\bsU}{\boldsymbol{U}}
\newcommand{\bsLam}{\boldsymbol{\Lambda}}
\newcommand{\bsu}{\boldsymbol{u}}
\newcommand{\bsv}{\boldsymbol{v}}
\newcommand{\bsw}{\boldsymbol{w}}
\newcommand{\bsf}{\boldsymbol{f}}
\newcommand{\bsg}{\boldsymbol{g}}
\newcommand{\bsh}{\boldsymbol{h}}

\newcommand{\bsr}{\boldsymbol{r}}

\title
[FFT based direct solver for the Helmholtz problem]
{A fast Fourier transform based direct solver for the Helmholtz problem}

\author{Jari Toivanen}
\author{Monika Wolfmayr}

\date{}

\thanks{This research was supported by the Academy of Finland, grant 295897.}

\address[J. Toivanen]{
Faculty of Information Technology; University of Jyväskylä;
P.O. Box 35 (Agora); FIN-40014 Jyväskylä; Finland}
\email{jari.a.toivanen@jyu.fi}

\address[M. Wolfmayr]{Faculty of Information Technology; University of Jyväskylä;
P.O. Box 35 (Agora); FIN-40014 Jyväskylä; Finland}
\email{monika.k.wolfmayr@jyu.fi}

\begin{document}

\maketitle

\begin{abstract}
This paper is devoted to the efficient numerical solution of
the Helmholtz equation in a two- or three-dimensional 
%orthogonal
rectangular
domain with an absorbing boundary condition (ABC).
The Helmholtz problem is discretized by standard bilinear and trilinear finite elements on an orthogonal mesh
yielding a separable system of linear equations.
The main key to high performance is to employ the Fast Fourier transform
(FFT) within a fast direct solver to solve the large separable systems.
The computational complexity of the proposed FFT based direct solver is $\mathcal{O}(N \log N)$ operations.
Numerical results for both two- and three-dimensional problems are presented
confirming the efficiency of the method discussed.
\end{abstract}

\section{Introduction}
\label{Sec:Intro}

This work presents an efficient fast direct solver 
employing FFT for the Helmholtz equation 
\begin{align}
\label{helmholtzIntro}
 - \triangle u - \omega^2 u = f 
\end{align}
in a rectangular domain with 
absorbing boundary conditions (ABCs).
The method can be applied for problems with the constant
zeroth order term coefficient $\omega$ being positive,
negative, or zero. Here, we focus on
large indefinite Helmholtz problems as they
are a common result from acoustic scattering
problems.
The zeroth order term coefficient
$\omega$ denotes the wave number, which is assumed to be a constant.
This condition is valid for homogeneous media.
For example, the iterative domain decomposition solution techniques
for acoustic scattering in layered media considered in \cite{HeiItoToi2018, Ito08} leads
a sequence of such problems.
In this work, we derive the fast solver for the case of having a first-order ABC imposed. 
However, the same method can be used for the second-order ABCs employed
in \cite{HeiRosToi2003}. These ABCs are the time-harmonic counterpart of the ABCs for
the wave equation considered in \cite{BamJolRob1990}. Moreover, the solver can be used
also for problems with Robin boundary conditions, which essentially
broadens its applicability.
Other separable Robin boundary conditions can be treated the same way as
the ABCs, which are just one possible type of boundary conditions for the Helmholtz equation
being solved by the method discussed in this work.
Solving the Helmholtz equation is in general difficult or impossible to solve efficiently with most numerical methods.
Difficulties related to the numerical solution of time-harmonic Helmholtz equations are discussed for instance
in \cite{Turkel2001}.

The numerical method proposed in this work is a fast direct solver,
which is applicable for problems posed in rectangular domains and
having suitable tensor product form matrices. 
This kind of diagonalization technique has already been proposed in
\cite{LynchRiceThomas1964}. 
However, we use FFT and sparsity in order to implement it efficiently.
The method is applicable for any
discretization leading to separable 9-point
stencil for two-dimensional problems and
27-point stencil for three-dimensional problems.
For example, bilinear
or trilinear finite element discretizations employed in this paper
lead to 
matrices with such suitable tensor product form.
Also the fourth-order accurate modified versions
of these elements \cite{Guddati04} are applicable with the solver.
When using fourth-order accurate bilinear or trilinear finite elements,
the reduced pollution error is obtained.
To improve the efficiency this method employs FFT instead of cyclic
reduction.
Several other fast direct solvers for elliptic problems
in rectangular domains are discussed, e.g., in \cite{Swarztrauber1977}.

The fast direct solvers can be used as
efficient preconditioners in the iterative solution of problems in more general domains,
see \cite{BuzDorGeoGol1971, HeiRosToi2003b}.
Other methods, which are based on an equivalent formulation
of the original problem to enable 
preconditioning with fast direct methods, are 
referred to as fictitious domain, domain imbedding, or capacitance matrix methods, 
and they have been successfully applied also to acoustic scattering problems 
for instance
in \cite{Ern1996, HeiKuzLip1999, HeiKuzNeiToi1998, HeiRosToi2003b}
as well as to scattering problems with multiple right-hand sides,
where a specific method for such problems is considered, for example,
in \cite{TezMacFar2001}.

A different 
fast direct solver for the Helmholtz equation 
\eqref{helmholtzIntro}
has already 
been presented in \cite{HeiRosToi2003}, 
where
cyclic reduction techniques are used in the solution procedure.
In \cite{HeiRosToi2003},
the application of a specific cyclic reduction-type fast direct method is presented, 
called the partial solution variant of the cyclic reduction (PSCR) method, 
to the solution of the Helmholtz equation 
\cite{RosToi1999b, RosToi1999, Vas1984}
and its computational efficiency is demonstrated.
For three-dimensional problems, a general fast direct solver like the 
PSCR 
has excellent parallel scalability,
see \cite{RosToi1999}.
However, the use of FFT over cyclic reduction techniques is preferable, 
since FFT is generally faster. That is partly due to very fast implementations of FFT.

The reason for not applying FFT previously was that 
the ABCs prevent the diagonalization using FFT.
The following three steps describe the basic idea to solve this problem
and are discussed in more detail in this work:
\begin{enumerate}
\item 
Some of the boundary conditions are changed on the ABC 
part of the boundary 
to be of periodic type. This
modified problem can be solved now with an FFT solver.
This step has the computational complexity of the FFT method, which is
$\mathcal{O}(N \log N)$ operations.
\item 
Since the boundary conditions have been changed to periodic ones,
the residual vector is computed due to these incorrect
boundary conditions.
This has nonzero components only on the boundary parts, where 
we have changed the boundary conditions.
The correction is computed by solving a problem with the original matrix,
where the right-hand side vector is 
the residual.
Only the component of the correction which lie on the changed boundaries are computed.
For this the so-called partial solution problem
a special technique exists, see 
\cite{Banegas78, HeiRosToi2003, Kuznetsov89, RosToi1999} for example.
Due to sparsity of the vectors the solution requires only $\mathcal{O}(N)$ or $\mathcal{O}(N \log N)$
operations in case of two or three dimensions, respectively,
where $N$ is the total number of unknowns.
After this step, the correct boundary values are known.
\item 
A similar problem to the one in the first step 
is solved, but now the right-hand side vector is adjusted so that the 
solution has correct values on the boundaries.
From this it follows that the solution is also the solution of the original problem.
The computational complexity is again $\mathcal{O}(N \log N)$ operations as in the first step.
\end{enumerate}

A similar method to the above one was already proposed in \cite{ErnGol1994},
but with the major difference that the problem in the second step 
was solved iteratively.
The new solver employing the FFT method is 
efficient in terms of computational time 
and memory usage,
especially in the three-dimensional case.
For three-dimensional problems, a general fast direct solver like the 
PSCR method
\cite{RosToi1999} requires $\mathcal{O}(N(\log N)^2)$ operations, whereas
the FFT based solver reduces the complexity to $\mathcal{O}(N \log N)$ operations.
The FFT method combined with the direct solver as proposed in this article
leads to the same, nearly optimal
complexity of $\mathcal{O}(N \log N)$ operations for both two- and three-dimensional problems.

The paper is organized as follows: In Section \ref{Sec:FormulationDiscretization}, we present the 
classic and variational formulation of the Helmholtz problem as well as its discretization by bilinear and trilinear finite elements
leading to a system of linear equations. The main idea of the solver for this problem and some preliminaries, which appear
in the initialization process of the implemented fast solver for the two- and three-dimensional case,
are discussed in Section \ref{Sec:MainIdeaPreliminaries}.
Sections \ref{Sec:FastSolver2D} and \ref{Sec:FastSolver3D} are devoted to the two- and three-dimensional problems, 
respectively. Here, some preliminaries as well as the fast solver steps are discussed.
Finally, numerical results for both two- and three-dimensional problems are presented in
Section \ref{Sec:Numerics}, and conclusions are drawn in Section \ref{Sec:Conclusions}.

\section{Problem formulation and discretization}
\label{Sec:FormulationDiscretization}

Let $\Omega \subset \mathbb{R}^d$, $d = 2,3$, be a $d$-dimensional rectangular domain and $\Gamma = \partial \Omega$
its boundary.
We consider the Helmholtz equation describing the linear propagation of time-harmonic acoustic waves given by
\begin{align}
 \label{helmholtz}
 - \triangle u - \omega^2 u &= f \qquad \text{in } \Omega, \\
 \label{bc}
 \mathcal{B} u &= 0 \qquad \text{on } \Gamma,
\end{align}
where $\omega$ denotes the wave number.
Here, equation \eqref{bc} is an approximation for the Sommerfeld radiation condition
\begin{align}
 \lim_{r \rightarrow \infty} r^{(d-1)/2} \left( \partial_r u - i \omega u \right) = 0
\end{align}
appearing for the general model problem in 
$\mathbb{R}^d$. The approximation is performed by 
truncating the unbounded domain 
at a finite distance, which provides 
the boundary condition at the truncation boundary.
In order to reduce spurious reflections caused by this artificial boundary, we use local absorbing boundary conditions.
In general, the boundary $\Gamma$ can be decomposed into parts where different boundary conditions are imposed.
Let now $\Gamma = \Gamma_B \cup \Gamma_N$ be decomposed into an absorbing boundary condition part
$\Gamma_B$ and a Neumann boundary condition part $\Gamma_N$.
The Neumann boundary conditions are given by 
\begin{align}
 \label{Neumannbc}
 \mathcal{B} u = \nabla u \cdot \bsn = 0 \qquad \text{on } \Gamma_N,
\end{align}
whereas in case of absorbing boundary conditions we have that 
\begin{align}
 \label{1st-order-abc}
 \mathcal{B} u = \nabla u \cdot \bsn - i \omega u = 0 \qquad \text{on } \Gamma_B,
\end{align}
where $\bsn$ denotes the outward normal to the boundary. 
This type of absorbing boundary conditions \eqref{1st-order-abc} are called
first-order absorbing boundary conditions.

\begin{remark}
In practical applications, so-called second-order 
absorbing boundary conditions are also important and were considered in 
\cite{HeiRosToi2003} as well.
However, the choice of the absorbing boundary conditions does not 
have any impact on the performance of the proposed solver. 
Thus, the same method can be used with given second-order 
absorbing boundary conditions, see \cite{HeiRosToi2003}.
\end{remark}

In order to obtain the discrete version of the Helmholtz problem
\eqref{helmholtz}--\eqref{bc}, let us first state its weak formulation.
For that we introduce the Hilbert space $V = H^1(\Omega)$. 
Let the source term $f \in L_2(\Omega)$ be given.  Multiplying \eqref{helmholtz} with a test function
as well as applying integration by parts and the boundary condition \eqref{bc}, yields the weak problem:
Find $u \in V$ such that
\begin{align}
 \label{helmholtzWeak}
 a(u,v) = \int_\Omega f v \, d\bsx \qquad \forall v \in V,
\end{align}
where
\begin{align}
 \label{bf-1st-order-abc}
a(u,v) = \int_\Omega \left( \nabla u \cdot \nabla v - \omega^2 u v \right) d\bsx
- i \omega \int_{\partial \Omega} u v \, ds.
\end{align}
Let us denote the mesh points by $x_{j,l}$, $l=1,\dots,n_j$ for every $x_j$-direction with
$j \in \{1,\dots,d\}$ and the corresponding mesh size
by $h_j = 1/(n_j-1)$.
Thus, the mesh is equidistant in each direction $x_j$.
Moreover, we denote the full mesh size by $N$, 
which is given by $N = n_1 \times \dots \times n_d$.
Discretizing problem \eqref{helmholtzWeak} by bilinear or trilinear finite elements
on an orthogonal mesh leads to a system of linear equations given by
\begin{align}
\label{probA}
\bsA \bsu = \bsf, 
\end{align}
where the matrix $\bsA$ has a separable tensor product form
and $\bsf = (f_1, \dots, f_N)$ denotes the discrete right-hand side.
For the two-dimensional case, the matrix $\bsA$ is given by 
\begin{align}
\label{A2D}
\bsA = (\bsK_1 - \omega^2 \bsM_1) \otimes \bsM_2 + \bsM_1 \otimes \bsK_2, 
\end{align}
whereas in three dimensions it is given by
\begin{align}
\label{A3D}
\bsA = (\bsK_1 - \omega^2 \bsM_1) \otimes \bsM_2 \otimes \bsM_3
+ \bsM_1 \otimes (\bsK_2 \otimes \bsM_3 + \bsM_2 \otimes \bsK_3). 
\end{align}
Here, the $n_j \times n_j$-matrices $\bsK_j$ and $\bsM_j$ are one-dimensional stiffness 
and mass matrices, respectively, in the $x_j$-direction 
with possible modifications on 
the boundaries due to the absorbing boundary conditions.
They are computed by one-dimensional numerical quadrature on the unit interval 
and are given as follows
\begin{align}
\label{matrixK}
\bsK_j = \frac{1}{h_j} \left( \begin{array}{rcccl}
k_{1,1}& -1 & & & \\
-1 & 2 & -1 & &  \\
 &  \ddots & \ddots & \ddots & \\
 & & -1 & 2 & -1 \\
 & & & -1 & k_{n_j,n_j} 
\end{array} \right)
\end{align}
and
\begin{align}
\label{matrixM}
\bsM_j = \frac{h_j}{6} \left( \begin{array}{ccccc}
2 & 1 & & & \\
1 & 4 & 1 & & \\
 &  \ddots & \ddots & \ddots & \\
 & & 1 & 4 & 1 \\
 & &  & 1 & 2
\end{array} \right)
\end{align}
where the first and last entries are including the corresponding boundary conditions.
Absorbing boundary conditions \eqref{1st-order-abc} yield the entries
$k_{1,1} = k_{n_j,n_j} = 1 - i \omega h_j$, whereas
Neumann boundary conditions lead to
$k_{1,1} = k_{n_j,n_j} = 1$. The matrix $\bsM_j$ is the same for both Neumann and (first-order)
absorbing boundary conditions.

\begin{remark}
\label{Rem:ABCx1}
Without loss of generality let us assume that the absorbing boundary conditions
are posed in direction of $x_1$ for both (opposite) sides.
\end{remark}

In the next two sections, an efficient fast solver for solving the large linear systems
$\bsA \bsu = \bsf$ is proposed. We split the discussion into two sections corresponding to
the two-dimensional and three-dimensional problems with the
matrix $\bsA$ given by \eqref{A2D} and \eqref{A3D}, respectively.

\section{Main idea and preliminaries}
\label{Sec:MainIdeaPreliminaries}

\subsection{Basic steps}

The idea for solving the problem $\bsA \bsu = \bsf$ 
is to consider an auxiliary problem 
\begin{align}
\label{probB}
\bsB \bsv = \bsf, 
\end{align}
where the system matrix 
$\bsB$ is derived by changing some 
absorbing boundary conditions
to periodic ones.
The key is that we can solve the modified (periodic) problem $\bsB \bsv = \bsf$ now 
by using the FFT method, which is not possible for the original problem $\bsA \bsu = \bsf$.

Before going into more details later, let us discuss the main steps of the solver.

\textit{Step 1:} 
Solve problem \eqref{probB} $\bsB \bsv = \bsf$. 

\textit{Step 2:} 
Let $\bsw = \bsu - \bsv$.
Solve 
\begin{align}
\label{probStep2}
\bsA \bsw = \bsf - \bsA \bsv = \bsB \bsv - \bsA \bsv = (\bsB - \bsA) \bsv.
\end{align}

\textit{Step 3:} 
Solve
\begin{align}
\label{probStep3}
\bsB \bsu  = \bsf + (\bsB - \bsA) (\bsv + \bsw).
\end{align}

Note that after Step 2, we would have already had 
the identity $\bsu = \bsv + \bsw$, but we do not use this
identity directly to compute $\bsu$, which will be explained 
in Subsections \ref{SSec:FS2d} and \ref{SSec:FS3d}.

\subsection{Preliminaries}

We have assumed that the absorbing boundary conditions are given in $x_1$-direction on both opposite boundaries,
see Remark \ref{Rem:ABCx1}.
Since we have to change the boundary conditions on the two opposite boundaries to be of periodic type
for the auxiliary problem \eqref{probB}, we consider the one-dimensional stiffness and mass matrices 
$\bsK_1$ and $\bsM_1$ and change the entries corresponding to the boundary parts into periodic conditions
for the auxiliary problem \eqref{probB}. We denote these matrices by
$\bsK_1^B$ and $\bsM_1^B$. The $n_j \times n_j$ circulant matrices $\bsK_j^B$ and $\bsM_j^B$
are given as follows
\begin{align}
\label{matrixKper}
\bsK_j^B = \frac{1}{h_j} \left( \begin{array}{ccccc}
2 & -1 & & & - 1 \\
-1 & 2 & -1 & &  \\
 &  \ddots & \ddots & \ddots & \\
 & & -1 & 2 & -1 \\
-1 & & & -1 & 2
\end{array} \right)
\end{align}
and
\begin{align}
\label{matrixMper}
\bsM_j^B = \frac{h_j}{6} \left( \begin{array}{ccccc}
4 & 1 & & & 1 \\
1 & 4 & 1 & & \\
 &  \ddots & \ddots & \ddots & \\
 & & 1 & 4 & 1 \\
1  & &  & 1 & 4
\end{array} \right)
\end{align}
which means that we have changed the boundary conditions on two opposite boundaries to be of periodic type.
Hence,
\begin{align}
\label{circulantvalues}
\begin{aligned}
\bsK^B_{j,(1,2)} &= \bsK^B_{j,(1,n_j)} = \bsK^B_{j,(2,1)} = \bsK^B_{j,(n_j,1)}, \\
\bsM_{j,(1,2)}^B &= \bsM_{j,(1,n_j)}^B = \bsM_{j,(2,1)}^B = \bsM_{j,(n_j,1)}^B.
\end{aligned}
\end{align}
Let us consider now the original problem \eqref{probA} and the auxiliary problem \eqref{probB} in more detail.
After a suitable permutation $\bsA$ and $\bsB$ have the block forms
\begin{align}
\bsA = \begin{pmatrix}
\bsA_{bb} & \bsA_{br} \\
\bsA_{rb} & \bsA_{rr}
\end{pmatrix}
\qquad\text{and}\qquad
\bsB = \begin{pmatrix}
\bsB_{bb} & \bsA_{br} \\
\bsA_{rb} & \bsA_{rr}
\end{pmatrix},
\end{align}
the subscripts $b$ and $r$ correspond to the nodes on the $\Gamma_B$ boundary
and to the rest of the nodes, respectively.
We denote by $|b|$ and $|r|$ the corresponding sizes of the set of nodes on the
$\Gamma_B$ boundary and of the rest of the nodes such that
$N = |b|+|r|$.
Note that the matrix $\bsB - \bsA$ has the structure
\begin{align}
\label{BminA}
\bsB - \bsA = \begin{pmatrix}
\bsB_{bb} - \bsA_{bb} & 0 \\
0 & 0
\end{pmatrix},
\end{align}
and only the matrix
\begin{align}
\label{Cbb}
\bsC_{bb} = \bsB_{bb} - \bsA_{bb}
\end{align}
has to be saved for the application of the fast solver.
Steps of the fast solver for $\bsA \bsu = \bsf$ include also the application of
the so-called partial solution method, 
which is a special implementation of the method of separation of variables.
This method involves the solution of the generalized eigenvalue problems
\begin{align}
\label{EVP1a}
\bsK_1 \bsV_1 = \bsM_1 \bsV_1 \bsLam^A_1 
\end{align}
and
\begin{align}
\label{EVP1b}
\bsK_1^B \bsW_1 = \bsM_1^B \bsW_1 \bsLam^B_1,
\end{align}
where the matrices $\bsLam^A_1$ and 
$\bsLam^B_1$ contain the eigenvalues 
as diagonal entries and the 
matrices $\bsV_1$ and $\bsW_1$ 
contain the corresponding eigenvectors as their columns.
Let us define the more general problem
\begin{align}
\label{EVPjb}
\bsK_j^B \bsW_j = \bsM_j^B \bsW_j \bsLam^B_j
\end{align}
with circulant matrices $\bsK_j^B$ and $\bsM_j^B$
for the problem size $n_j$ corresponding to any $x_j$-direction.
The $n_j \times n_j$ eigenvector matrix $\bsW_j$ 
for the generalized eigenvalue problem 
\eqref{EVPjb} is given by
\begin{align}
\label{eigvectorsPer}
\bsW_j = \left( \begin{array}{cccccc}
1 & 1 & 1 & \dots & 1 & 1 \\
1 & e^{-\frac{2 \pi i}{n_j}} & e^{-2 \frac{2 \pi i}{n_j}} & \dots & e^{-(n_j-2) \frac{2 \pi i}{n_j}} &  e^{-(n_j-1) \frac{2 \pi i}{n_j}} \\
1 & e^{-2 \frac{2 \pi i}{n_j}} & e^{-4 \frac{2 \pi i}{n_j}} & \dots & e^{-2 (n_j-2) \frac{2 \pi i}{n_j}} &  e^{-2 (n_j-1) \frac{2 \pi i}{n_j}} \\
 \vdots & \vdots  & \vdots & \ddots & \vdots &  \vdots \\
1 & e^{-(n_j-2) \frac{2 \pi i}{n_j}} & e^{-(n_j-2) 2 \frac{2 \pi i}{n_j}} & \dots & e^{-(n_j-2)^2 \frac{2 \pi i}{n_j}} &  e^{-(n_j-1) (n_j-2) \frac{2 \pi i}{n_j}} \\
1 & e^{-(n_j-1) \frac{2 \pi i}{n_j}} & e^{-2(n_j-1) \frac{2 \pi i}{n_j}} & \dots & e^{-(n_j-2) (n_j-1) \frac{2 \pi i}{n_j}} &  e^{-(n_j-1)^2 \frac{2 \pi i}{n_j}} 
\end{array} \right)
\end{align}
and the diagonal entries of the corresponding diagonal eigenvalue matrix $\bsLam^B_j$ are 
\begin{align}
\label{eigvaluesPer}
\bsLam^B_{j,l} = \frac{\bsK_{j,(1,1)}^B + \bsK_{j,(1,n_j)}^B e^{-(l-1)\frac{2 \pi i}{n_j}} + \bsK_{j,(1,2)}^B e^{-(l-1)(n_j-1) \frac{2 \pi i}{n_j}} }{
\bsM_{j,(1,1)}^B + \bsM_{j,(1,n_j)}^B e^{-(l-1)\frac{2 \pi i}{n_j}} + \bsM_{j,(1,2)}^B e^{-(l-1)(n_j-1) \frac{2 \pi i}{n_j}} }
\end{align}
for $l = 1,\dots,n_j$.
In order to apply the partial solution method the eigenvectors are
normalized so that they satisfy the conditions:
\begin{align}
\label{EVPcondA}
\bsV_1^T \bsM_1  \bsV_1 = \boldsymbol{I}_1
\qquad &\text{ and } \qquad
\bsV_1^T \bsK_1  \bsV_1 = \bsLam^A_1, \\
\label{EVPcondB}
\bsW_1^T \bsM_1^B  \bsW_1 = \boldsymbol{I}_{1}
\qquad &\text{ and } \qquad
\bsW_1^T \bsK_1^B  \bsW_1 = \bsLam^B_1. 
\end{align}
This can be achieved by multiplying $\bsV_1$ and
$\bsW_1$ with the vectors $\bss^A_1$ and $\bss^B_1$ of length $n_1$, respectively.
We denote the general vectors by $\bss^A_j$ and $\bss^B_j$ of length $n_j$.
The entries of $\bss^A_j$ are given by 
\begin{align}
\label{sA}
\bss_{j,l}^A = 1/\|\bsM_j\|_{\bsV_{j,l}}
\end{align}
with the componentwise matrix norms
\begin{align}
\label{sAnorm}
\|\bsM_j\|_{\bsV_{j,l}} = \left(\bsV_{j,l}^T \bsM_j  \bsV_{j,l} \right)^{1/2},
\end{align}
which are weighted by the $l$-th eigenvectors of $\bsV_j$ for $l = 1,\dots,n_j$.
Similarly, the entries of $\bss^B_j$ are given by 
\begin{align}
\label{sB}
\bss_{j,l}^B = 1/\|\bsM_j^B\|_{\bsW_{j,l}}
\end{align}
with the componentwise matrix norms
\begin{align}
\label{sBnorm}
\|\bsM_j^B\|_{\bsW_{j,l}} = \left(\bsW_{j,l}^T \bsM_j^B  \bsW_{j,l} \right)^{1/2},
\end{align}
which are weighted by the $l$-th eigenvectors of $\bsW_j$ for $l = 1,\dots,n_j$.
In \eqref{EVPcondA} and \eqref{EVPcondB}, 
$\boldsymbol{I}_{1}$ denotes the identity 
matrix of size $n_1 \times n_1$. 
In the following, $\boldsymbol{I}_{j}$ and $\boldsymbol{I}_{jk}$
denote the identity matrices of 
sizes 
$n_j \times n_j$ and $(n_j \times n_k) \times (n_j \times n_k)$, respectively.
The eigenvector matrices $\bsV_1$ and 
$\bsW_1$ are used for solving the partial solution problems 
when needed in the steps of the solver.
However, the generalized eigenvalue problems \eqref{EVP1a} and \eqref{EVP1b}
have to be solved only once during the solution process -- in the initialization.
The conditions \eqref{EVPcondA} and 
\eqref{EVPcondB}
also lead to a convenient representation for the inverses of the system matrices
$\bsA$ and $\bsB$, which is discussed in Subsections
\ref{SSec:Ini2d} and \ref{SSec:Ini3d} for the respective two-dimensional 
and three-dimensional case.
In the next two sections, we discuss the two-dimensional 
and three-dimensional 
problems in more detail separately,
since the efficient implementation of the initialization process and the steps of the fast solver 
differs in both cases.

\section{The two-dimensional case} 
\label{Sec:FastSolver2D}

\subsection{Preliminaries for the two-dimensional problem}
\label{SSec:Ini2d}

\subsubsection{Reformulation of problem matrices $\bsA$ and $\bsB$}

The matrix $\bsB$ for the auxiliary problem \eqref{probB} in the 
two-dimensional case is given by
\begin{align}
\label{B2D}
\bsB = (\bsK_1^B - \omega^2 \bsM_1^B) \otimes \bsM_2 + \bsM_1^B \otimes \bsK_2. 
\end{align}
Using the conditions 
\eqref{EVPcondB} as follows
\begin{align}
\label{B2DinvProcess}
\begin{aligned}
\bsB &= (\bsW_1^{-T} \bsLam_1^B \bsW_1^{-1} - \omega^2 \bsW_1^{-T} \boldsymbol{I}_{1} \bsW_1^{-1}) \otimes \bsM_2 
+ (\bsW_1^{-T} \boldsymbol{I}_{1} \bsW_1^{-1}) \otimes \bsK_2 \\
&= (\bsW_1^{-T} (\bsLam_1^B  - \omega^2  \boldsymbol{I}_{1} ) \bsW_1^{-1}) \otimes \bsM_2 
+ (\bsW_1^{-T} \boldsymbol{I}_{1} \bsW_1^{-1}) \otimes \bsK_2 \\
&= (\bsW_1^{-T} \otimes \bsI_{2}) ((\bsLam_1^B  - \omega^2  \boldsymbol{I}_{1}) \otimes \bsM_2 
+ \boldsymbol{I}_{1} \otimes \bsK_2) (\bsW_1^{-1} \otimes \bsI_{2}),
\end{aligned}
\end{align}
the inverse of $\bsB$ can be represented by
\begin{align}
\label{B2Dinv}
\bsB^{-1} = (\bsW_1 \otimes \bsI_{2}) \, \bsH_B^{-1}
(\bsW_1^T \otimes \bsI_{2}),
\end{align}
where
\begin{align}
\label{HB2D}
\bsH_B = (\bsLam_1^B - \omega^2 \bsI_{1}) \otimes \bsM_2 + \bsI_{1} \otimes \bsK_2.
\end{align}
Similarly using the conditions \eqref{EVPcondA}, we obtain for the system matrix $\bsA$ the following representation:
\begin{align}
\label{A2Dinv}
\begin{aligned}
\bsA ^{-1}
&= (\bsV_1 \otimes \bsI_2) \, \bsH_A^{-1} (\bsV_1^T \otimes \bsI_2), 
\end{aligned}
\end{align}
where
\begin{align}
\label{HA2D}
\bsH_A = (\bsLam_1^A - \omega^2 \bsI_{1}) \otimes \bsM_2 + \bsI_{1} \otimes \bsK_2.
\end{align}

\subsubsection{LU decomposition}
\label{SSec:LU2d}

Linear systems with the block diagonal matrices $\bsH_A$ and $\bsH_B$ are
solved with a direct method.
The matrices $\bsH_A$ and $\bsH_B$ consist of $N$ 
diagonal blocks each of size $N$.
In the following, let us discuss 
the method applied on the matrix $\bsH_B$ (given by \eqref{HB2D}),
since it applies analogously for $\bsH_A$.
First,  the LU decomposition of $\bsH_B$ 
is computed as follows
\begin{align}
\label{HBLU2d}
\bsH_B = \bsL_B \bsU_B.
\end{align}
Then the linear system
\begin{align}
\label{HB2DLU1}
\bsH_B \bsy = \bsL_B \bsU_B \bsy = \bsr
\end{align}
is solved by
solving
the respective two subproblems 
\begin{align}
\label{HB2DLU2}
\bsL_B \bsz = \bsr \qquad \text{ and } \qquad 
\bsU_B \bsy = \bsz
\end{align}
consecutively in the
application of the fast solver.
Note that $\bsr$ denotes now some right-hand side which is in the different steps of the solver also different.
However, we will describe that in more detail in the next subsection. 
We denote by
\begin{align}
\label{HALU2d}
\bsH_A = \bsL_A \bsU_A
\end{align}
the LU decomposition corresponding to $\bsH_A$.
The structure of $\bsH_B$ and $\bsH_A$ is essential for the fast application of the direct solver.
The diagonal blocks of these matrices are tridiagonal which makes the LU decomposition
fast for them.
The computational complexity is optimal
$\mathcal{O}(N)$. 

\subsection{Fast solver in the two-dimensional case}
\label{SSec:FS2d}

\subsubsection{Step 1}

The auxiliary problem
\begin{align}
\label{step1}
\bsB \bsv = \bsB \begin{pmatrix} \bsv_b \\ \bsv_r \end{pmatrix} = \bsf
\end{align}
is solved, but only $\bsv_b$ and not $\bsv_r$ is computed.
For that, the FFT $\hat \bsf$ of the 
right-hand side $\bsf$ is computed first, 
where its coefficients $\hat f_k$ are given as follows
\begin{align}
\label{step1FFTfcoeff}
\hat f_k = \sum_{l = 1}^N e^{-2 \pi i \frac{(l-1) (k-1)}{N}} f_l
\end{align}
for all $k = 1, \dots, N$ and 
normalization might be needed, 
which means that $\hat \bsf$
is multiplied by the vector $\bss^B_1$, where its components are defined in \eqref{sB} with $j=1$.
Note that computing the FFT corresponds to the multiplication $\hat \bsf = (\bsW_1^T \otimes \bsI_{2}) \, \bsf$.
The vector $\hat \bsf$ is saved, since it will be needed in Step 3
as well. 
Next, we apply the LU decomposition \eqref{HB2DLU1} with the right-hand side
$\hat \bsf$ 
\begin{align}
\label{2dLUBs1}
\bsL_B \bsz_1 = \hat \bsf \qquad \text{ and } \qquad 
\bsU_B \tilde \bsz_1 =  \bsz_1.
\end{align}
Performing the inverse FFT 
on the resulting vector $\tilde \bsz_1$
would provide both $\bsv_b$ and $\bsv_r$, which would
correspond to the multiplication $\bsv = (\bsW_1 \otimes \bsI_{2}) \, \tilde \bsz_1$.
However, since we only need $\bsv_b$, instead of that,
we multiply the vector $\tilde \bsz_1$ 
by the matrix $\bsW_1$ 
by taking advantage of the sparsity of the desired components $\bsv_b$.
More precisely, we multiply $\tilde \bsz_1$ 
from the left side by the eigenvectors of $\bsW_1$ which 
correspond only to the boundary
$\Gamma_B$ denoted by 
the matrix $\bsW_1^b$
of size $|b| \times n_1$
leading to 
\begin{align}
\label{2dBs1}
\bsv_b = (\bsW_1^b \otimes \bsI_{2}) \, \tilde \bsz_1,
\end{align}
which resembles the representation \eqref{B2Dinv} for $\bsB^{-1}$.
The computational complexity of Step 1 is $\mathcal{O}(N \log N)$.

\subsubsection{Step 2}

We introduce the additional vector $\bsw$ defined as 
$\bsw = \bsu - \bsv$.
Since
\begin{align}
\label{step2decomp}
\begin{aligned}
\bsA \bsw = \bsA \bsu - \bsA \bsv = \bsf - \bsA \bsv = \bsB \bsv - \bsA \bsv
\end{aligned}
\end{align}
and \eqref{BminA},
we obtain the following problem:
\begin{align}
\label{step2}
\begin{aligned}
\bsA \bsw = \bsA \begin{pmatrix} \bsw_b \\ \bsw_r \end{pmatrix} 
= (\bsB - \bsA) \bsv
= \begin{pmatrix} \bsC_{bb} 
 \bsv_b \\ 0 \end{pmatrix},
\end{aligned}
\end{align}
where $\bsC_{bb}$ is defined by \eqref{Cbb}.
We solve the problem \eqref{step2} by using the representation \eqref{A2Dinv},
but compute again only $\bsw_b$ and not $\bsw_r$, since
the right-hand side as well as the desired components of the solution are both sparse.
First, we compute
\begin{align}
\label{step2-1}
\begin{aligned}
 \bsg_b = ({\bsV_1^b}^T \otimes \bsI_2) \, \bsC_{bb}  \bsv_b,
\end{aligned}
\end{align}
where 
the matrix $\bsV_1^b$
of size $|b| \times n_1$
denotes the eigenvectors of $\bsV_1$ which 
correspond only to the boundary
$\Gamma_B$.
Then
by using the LU decomposition \eqref{HALU2d} we solve
\begin{align}
\label{step2-2}
\begin{aligned}
\bsL_A \bsz_2 = \bsg_b \qquad \text{ and } \qquad 
\bsU_A \tilde \bsz_2 =  \bsz_2,
\end{aligned}
\end{align}
and, finally, we solve
\begin{align}
\label{step2-3}
\begin{aligned}
 \bsw_b = (\bsV_1^b \otimes \bsI_2) \, \tilde \bsz_2,
\end{aligned}
\end{align}
which resembles the representation \eqref{A2Dinv} for $\bsA^{-1}$ 
but corresponds only to
the boundary $\Gamma_B$. 
The computational complexity of Step 2 is of optimal order
$\mathcal{O}(N)$.

\subsubsection{Step 3}

Finally, in order to obtain the solution $\bsu$ of the original problem \eqref{probA},
we solve now the problem 
\begin{align}
\label{step3}
\begin{aligned}
\bsB \bsu = \bsf + (\bsB - \bsA) (\bsv + \bsw)
= \bsf + \begin{pmatrix} \bsC_{bb} 
(\bsv_b + \bsw_b) \\ 0 \end{pmatrix}
\end{aligned}
\end{align}
due to 
$\bsB \bsu = \bsA \bsu + \bsB \bsu - \bsA \bsu$.
Since we have already computed 
the Fourier transformation $\hat \bsf$ of $\bsf$ in Step 1 by \eqref{step1FFTfcoeff},
we only need to compute 
the Fourier transformation 
of the second term of the right-hand side of equation \eqref{step3} and due to sparsity again only 
the part corresponding to the boundary $\Gamma_B$ as follows
\begin{align}
\label{step3term}
\begin{aligned}
\bsh_b = ({\bsW_1^b}^T \otimes \bsI_2) \, \bsC_{bb} (\bsv_b + \bsw_b)
\end{aligned}
\end{align}
leading to the Fourier transformation 
of the entire right-hand side of equation \eqref{step3}
denoted by $\hat \bsf + \hat \bsh$.
Next, we apply the LU decomposition \eqref{HB2DLU1} with this right-hand side yielding
\begin{align}
\label{2dLUBs3}
\bsL_B \bsz_3 = \hat \bsf + \hat \bsh \qquad \text{ and } \qquad 
\bsU_B \tilde \bsz_3 =  \bsz_3.
\end{align}
In the last step all resulting components 
are needed (not only the ones corresponding to $\Gamma_B$) to obtain the solution $\bsu$
by applying the inverse Fourier transformation 
on $\tilde \bsz_3$,
where its coefficients $\hat u_k$ are given as follows
\begin{align}
\label{step3IFFTz3coeff}
\hat u_k = \frac{1}{N} \sum_{l = 1}^N e^{2 \pi i \frac{(l-1) (k-1)}{N}} \tilde z_{3,l}
\end{align}
for all $k = 1, \dots, N$. Again 
the vector $\tilde \bsz_3$ may need to be
multiplied by the vector $\bss^B_1$ with entries 
defined in 
\eqref{sB} for $n_1$ 
before applying the inverse Fourier transformation on it.
The computation of the inverse Fourier transformation 
corresponds to the multiplication 
\begin{align}
\label{2dBs3}
\bsu = (\bsW_1 \otimes \bsI_{2}) \, \tilde \bsz_3.
\end{align}
The computational complexity of Step 3 is $\mathcal{O}(N \log N)$.

\begin{conclusion}
\label{conclusion2d}
Combing the results on the computational complexities of all three steps finally leads to the overall computational complexity
of $\mathcal{O}(N \log N)$ operations for the two-dimensional problem.
\end{conclusion}

\section{The three-dimensional case} 
\label{Sec:FastSolver3D}

\subsection{Preliminaries for the three-dimensional problem}
\label{SSec:Ini3d}

\subsubsection{Matrices $\bsA$ and $\bsB$}

The matrix $\bsB$ in the three-dimen\-sional 
case is given by
\begin{align}
\label{B3D}
\bsB = (\bsK_1^B - \omega^2 \bsM_1^B) \otimes \bsM_2 \otimes \bsM_3
+ \bsM_1^B \otimes (\bsK_2 \otimes \bsM_3 + \bsM_2 \otimes \bsK_3).
\end{align}
Using the equations \eqref{EVPcondA} and 
\eqref{EVPcondB},
the inverses of $\bsA$ and $\bsB$ can be represented 
analogously as in the two-dimensional case as follows
\begin{align}
\label{A3Dinv}
\bsA^{-1}
= (\bsV_1 \otimes \bsI_{23}) \, \bsH_A^{-1} (\bsV_1^T \otimes \bsI_{23}) 
\end{align}
and
\begin{align}
\label{B3Dinv}
\bsB^{-1} = (\bsW_1 \otimes \bsI_{23}) \, \bsH^{-1}_B
(\bsW_1^T \otimes \bsI_{23}), 
\end{align}
where
\begin{align}
\label{HA3D}
\bsH_A = ((\bsLam^A_1 - \omega^2 \bsI_{1}) \otimes \bsM_2 + \bsI_{1} \otimes \bsK_2) \otimes \bsM_3
+ \bsI_{1} \otimes \bsM_2 \otimes \bsK_3 
\end{align}
and
\begin{align}
\label{HB3D}
\bsH_B = ((\bsLam^B_1 - \omega^2 \bsI_{1}) \otimes \bsM_2 + \bsI_{1} \otimes \bsK_2) \otimes \bsM_3
+ \bsI_{1} \otimes \bsM_2 \otimes \bsK_3, 
\end{align}
respectively.

\subsubsection{Applying the two-dimensional fast solver}
\label{SSec:HAHB3d}

Linear systems with the block diagonal matrices $\bsH_B$ and $\bsH_A$ are again
solved with a direct method.
However, 
it is performed in a different way than for the two-dimensional 
case,
since the LU decomposition is 
slow for large block tridiagonal problems.
More precisely, the efficient implementation 
in three dimensions 
contains the application of 
the two-dimensional 
fast solver for $n_1$ subproblems of size $n_2 \times n_3$ including
the computation of the partial solution method 
in $x_2$-direction. 
For that, one needs to solve
the following generalized eigenvalue 
problems during the initialization process:
\begin{align}
\label{EVP2}
\bsK_2 \bsV_2 = \bsM_2 \bsV_2 \bsLam^A_2 \qquad \text{ and } \qquad
\bsK_2^B \bsW_2 = \bsM_2^B \bsW_2 \bsLam^B_2
\end{align}
with the matrices
$\bsK_2^B$ and $\bsM_2^B$ 
as defined in \eqref{matrixKper} and \eqref{matrixMper} but for the $x_2$-direction now.
Moreover, the diagonal eigenvalue matrices
$\bsLam^A_2$ and $\bsLam^B_2$ 
and the eigenvector matrices $\bsV_2$ and $\bsW_2$
are formed analogously as for the $x_1$-direction only the problem size changes to $n_2$.
For that, the generalized eigenvalue problem with circulant matrices for a general problem size $n_j$ was defined in
\eqref{EVPjb} and the corresponding eigenvector matrix is represented in \eqref{eigvectorsPer}. 
The eigenvectors are normalized to satisfy the conditions
\begin{align}
\label{EVPcondA2}
\bsV_2^T \bsM_2  \bsV_2 = \boldsymbol{I}_2
\qquad &\text{ and } \qquad
\bsV_2^T \bsK_2  \bsV_2 = \bsLam^A_2, \\
\label{EVPcondB2}
\bsW_2^T \bsM_2^B  \bsW_2 = \boldsymbol{I}_2
\qquad &\text{ and } \qquad
\bsW_2^T \bsK_2^B  \bsW_2 = \bsLam^B_2. 
\end{align}
This can be achieved by multiplying $\bsV_2$ and
$\bsW_2$ with the vectors $\bss^A_2$ and $\bss^B_2$ of length $n_2$, respectively,
which are 
 introduced in \eqref{sA}--\eqref{sBnorm}.
The fast solver includes the application of the partial solution method $n_1$ times
for the following two-dimensional 
$n_2 \times n_3$ subproblems in $x_2$-direction
corresponding to the matrices \eqref{HA3D} and \eqref{HB3D}:
\begin{align}
\label{A2D2}
\bsA_{A,l} = (\bsK_2 - \underbrace{(\omega^2 - \bsLam^A_{1,l})}_{=:p_A} \bsM_2) \otimes \bsM_3 + \bsM_2 \otimes \bsK_3
\end{align}
and
\begin{align}
\label{B2D2}
\bsA_{B,l} = (\bsK_2 - \underbrace{(\omega^2 - \bsLam^B_{1,l})}_{=:p_B} \bsM_2) \otimes \bsM_3 + \bsM_2 \otimes \bsK_3,
\end{align}
where $l=1,\dots,n_1$.
This yields the direct 
solution for $\bsH_A^{-1}$ and $\bsH_B^{-1}$.
Note that the subproblems \eqref{A2D2} 
and \eqref{B2D2}
reflect the structure of \eqref{HA3D} and \eqref{HB3D}, respectively, 
in the framework of 
the two-dimensional 
problem \eqref{A2D}. 
Now the parameters $p_A$ and $p_B$ have been introduced
in \eqref{A2D2} and \eqref{B2D2}, respectively.
The corresponding auxiliary problems are given by
\begin{align}
\label{A2D2aux}
\bsB_{A,l} = (\bsK_2^B - p_A \bsM_2^B) \otimes \bsM_3 + \bsM_2^B \otimes \bsK_3,
\end{align}
and
\begin{align}
\label{B2D2aux}
\bsB_{B,l} = (\bsK_2^B - p_B \bsM_2^B) \otimes \bsM_3 + \bsM_2^B \otimes \bsK_3,
\end{align}
which now reflect the structure of \eqref{HA3D} and \eqref{HB3D}, respectively,
in the framework of 
the two-dimensional 
problem \eqref{B2D}.
As defined in \eqref{eigvaluesPer},
$\bsLam^B_{1,l}$ is the $l$th-diagonal entry of the diagonal eigenvalue matrix
$\bsLam^B_{1}$. Analogously, we have denoted by 
$\bsLam^A_{1,l}$ the $l$th-diagonal entry of the diagonal eigenvalue matrix
$\bsLam^A_{1}$.

In summary, 
the two-dimensional 
problems \eqref{A2D2} and  \eqref{B2D2} are solved for all $l=1,\dots,n_1$
by applying  $n_1$ times the two-dimensional 
FFT based fast solver from Subsection \ref{SSec:FS2d} using 
the auxiliary two-dimensional 
problems \eqref{A2D2aux} and \eqref{B2D2aux}.
The computational complexity of this step is now 
$\mathcal{O}(N \log N)$, %.
since the FFT method has to be applied now in this step as well.

\subsection{Fast solver in the three-dimensional 
case}
\label{SSec:FS3d}

The three-dimensional 
solver has the same three main steps, but now instead of applying LU decomposition we apply the two-dimensional 
solver as described in the previous subsection.
We present all the main steps of 
the solver similarly as for the two-dimensional case in Subsection \ref{SSec:FS2d}
but in less detail. Hence, we refer the reader also to Subsection \ref{SSec:FS2d} for more details.

\subsubsection{Step 1}

The auxiliary problem \eqref{step1}
\begin{align*}
\bsB \bsv = \bsB \begin{pmatrix} \bsv_b \\ \bsv_r \end{pmatrix} = \bsf
\end{align*}
is solved, but only $\bsv_b$ and not $\bsv_r$ is computed.
For that, we need to compute $\hat \bsf = (\bsW_1^T \otimes \bsI_{23}) \, \bsf$ first,
which is equivalent to the computation of the FFT $\hat \bsf$  
for the right-hand side $\bsf$, 
where its coefficients $\hat f_k$ are given in \eqref{step1FFTfcoeff}
for all $k = 1, \dots, N$. 
The vector $\hat \bsf$ can be normalized by
multiplying it with the vector $\bss^B_1$ of length $n_1$
with components defined in \eqref{sB}.
The vector $\hat \bsf$ is saved, since it will be needed in Step 3 as well. 
Next, we apply $n_1$ times the two-dimensional 
solver for the $n_2 \times n_3$ problems
\eqref{B2D2} with the auxiliary problems \eqref{B2D2aux} leading to the solution
of the problem
\begin{align}
\label{3dHBs1}
\bsH_B \tilde \bsz_1 = \hat \bsf.
\end{align}
Since we only need $\bsv_b$, 
we now multiply the vector $\tilde \bsz_1$ 
by the matrix $\bsW_1$ 
by taking advantage of the sparsity of the desired components $\bsv_b$.
More precisely, we multiply $\tilde \bsz_1$ 
from the left side by 
the components of the eigenvectors of $\bsW_1$ which 
correspond only to the boundary
$\Gamma_B$ denoted by $\bsW_1^b$
leading to 
\begin{align}
\label{3dBs1}
\bsv_b = (\bsW_1^b \otimes \bsI_{23}) \, \tilde \bsz_1,
\end{align}
which resembles the representation \eqref{B3Dinv} for $\bsB^{-1}$.

\subsubsection{Step 2}

We introduce the additional vector $\bsw$ defined as 
$\bsw = \bsu - \bsv$.
Since \eqref{step2decomp}
and \eqref{BminA},
we derive at problem \eqref{step2}
\begin{align*}
\begin{aligned}
\bsA \bsw = \bsA \begin{pmatrix} \bsw_b \\ \bsw_r \end{pmatrix} 
= (\bsB - \bsA) \bsv
= \begin{pmatrix} \bsC_{bb} 
 \bsv_b \\ 0 \end{pmatrix}.
\end{aligned}
\end{align*}
We solve this problem 
by using the representation \eqref{A3Dinv},
but compute again only $\bsw_b$ and not $\bsw_r$, since
the right-hand side as well as the desired components of the solution are both sparse.
First, we solve
\begin{align}
\label{step2-1-3d}
\begin{aligned}
 \bsg_b = ({\bsV_1^b}^T \otimes \bsI_{23}) \, \bsC_{bb}  \bsv_b,
\end{aligned}
\end{align}
then by applying $n_1$ times the two-dimensional 
solver for the $n_2 \times n_3$ subproblems
\eqref{A2D2} with the auxiliary problems \eqref{A2D2aux} leading to the solution
of the problem
\begin{align}
\label{step2-2-3d}
\bsH_A \tilde \bsz_2 = \bsg_b.
\end{align}
Finally, we compute 
\begin{align}
\label{step2-3-3d}
\begin{aligned}
 \bsw_b = (\bsV_1^b \otimes \bsI_{23}) \, \tilde \bsz_2,
\end{aligned}
\end{align}
which resembles the representation \eqref{A3Dinv} for $\bsA^{-1}$ 
but corresponds only to
the boundary $\Gamma_B$. 

\subsubsection{Step 3}

The last step yields 
the solution $\bsu$ of the original problem \eqref{probA}.
First, we solve the problem \eqref{step3}
\begin{align*}
\begin{aligned}
\bsB \bsu = \bsf + (\bsB - \bsA) (\bsv + \bsw)
= \bsf + \begin{pmatrix} \bsC_{bb} 
(\bsv_b + \bsw_b) \\ 0 \end{pmatrix}
\end{aligned}
\end{align*}
due to $\bsB \bsu = \bsA \bsu + \bsB \bsu - \bsA \bsu$.
Since we have already computed 
the Fourier transformation $\hat \bsf$ of $\bsf$ in Step 1 by \eqref{step1FFTfcoeff},
we only need to compute 
the Fourier transformation 
of the second term of the right-hand side of equation \eqref{step3} and due to sparsity again only 
the part corresponding to the boundary $\Gamma_B$ as follows
\begin{align}
\label{step3term3d}
\begin{aligned}
\bsh_b = ({\bsW_1^b}^T \otimes \bsI_{23}) \, \bsC_{bb} (\bsv_b + \bsw_b)
\end{aligned}
\end{align}
leading to the Fourier transformation 
of the entire right-hand side of equation \eqref{step3}
denoted by $\hat \bsf + \hat \bsh$.
Next, we apply $n_1$ times the two-dimensional 
solver for the $n_2 \times n_3$ problems
\eqref{B2D2} with the auxiliary problems \eqref{B2D2aux} leading to the solution
of the problem
\begin{align}
\label{3dHBs3}
\bsH_B \tilde \bsz_3 = \hat \bsf + \hat \bsh.
\end{align}
In the last step all resulting components 
are needed (not only the ones corresponding to $\Gamma_B$) to obtain the solution $\bsu$
by applying the inverse Fourier transformation on $\tilde \bsz_3$,
where its coefficients $\hat u_k$ are defined as in \eqref{step3IFFTz3coeff}
for all $k = 1, \dots, N$. 
The vector $\tilde \bsz_3$ may be normalized by
multiplying it with the vector $\bss^B_1$ of length $n_1$
with the components defined in \eqref{sB}
before applying the inverse Fourier transformation
which corresponds
to the multiplication 
\begin{align}
\label{3dBs3}
\bsu = (\bsW_1 \otimes \bsI_{23}) \, \tilde \bsz_3.
\end{align}
finally leading to the solution $\bsu$.

\begin{conclusion}
\label{conclusion3d}
The computational complexities in all three steps for the three-dimensional problem are
$\mathcal{O}(N \log N)$ operations leading again
to %an 
the overall complexity of $\mathcal{O}(N \log N)$ operations.
\end{conclusion}

In the next section, we present numerical experiments for both the two-dimensional 
and three-dimensional case.

\section{Numerical results}
\label{Sec:Numerics}

For $d = 2$ and $d=3$,
the computational domain is chose as 
$\Omega = [0,1]^d$
the unit square and unit cube, respectively.
The wave number is set 
$\omega = 2 \pi$ for all numerical experiments.
The discretization meshes are uniform with respect to
each $x_j$-direction,  
where the corresponding step sizes are denoted by $h_j = 1/(n_j-1), j = 1, \dots, d$.
The right-hand side is chosen as 0.01 for the first $n_1$ entries and 1 for all the other entries.
In this set of experiments, the efficiency of the fast direct solver is discussed. The numerical experiments
have been computed 
using
MATLAB 
9.3, R2017b.

In the following, we compare the CPU times in seconds for computing the solution by applying Matlab's backslash
and the fast solver
presented in this work.
In our case,
Matlab's backslash uses the sparse direct solver UMFPACK
for computing the solution of the sparse linear systems
\eqref{A2D} and \eqref{A3D},
see \cite{Davis2006}.
Besides the computational times 
of applying Matlab's backslash and the FFT based direct solver
also the computational times needed in the initialization process are presented.
However, the initialization
processes differ between the two- 
and three-dimensional 
problems, which has already been addressed in Sections
\ref{Sec:FastSolver2D} and
\ref{Sec:FastSolver3D}.
Summing up, in the two-dimensional case 
the LU decompositions are computed for the matrices $\bsH_B$ and
$\bsH_A$ defined in \eqref{HB2D} and \eqref{HA2D}, respectively, 
(see Subsection \ref{SSec:LU2d}),
whereas the three-dimensional 
problem is not solved using LU decomposition in three dimensions. 
In this case, the action of the inverses of $\bsH_B$ and
$\bsH_A$ defined in \eqref{HB3D} and \eqref{HA3D}, respectively, 
are computed by applying the fast solver in two dimensions $n_1$ times
(see Subsection \ref{SSec:HAHB3d}).
This approach leads to a more efficient implementation in three dimensions. 

The largest numerical experiments were computed in three dimensions 
with 135 005 697 $(= 513^3)$ unknowns.
Note that this size is already too large regarding the computational memory in Matlab
for setting up the whole matrix matrix $\bsA$ as defined in \eqref{A3D},
which would be needed in order to apply Matlab's backslash. However, 
the fast solver presented in this work still computes the solution in a reasonable amount of time
without using too much computational memory, since the fast solver does not need 
to form the whole matrix $\bsA$, but only of the left upper block \eqref{Cbb} denoted by
$\bsC_{bb}$ of the
matrix $\bsB - \bsA$ defined in
\eqref{BminA}.

\begin{remark}
\label{remarkLaptop}
All numerical experiments presented in this section
were performed on a laptop with Intel(R) Core(TM) i5-6267U CPU @ 2.90GHz processor
and 16 GB 2133 MHz LPDDR3 memory.
\end{remark}

\subsection{Numerical results for the 
two-dimensional 
case}
\label{Ssec:Num2d}

For the first set of numerical experiments, we choose $n = n_1 = n_2$. 
The CPU times in seconds for the initialization process as well as the application of the fast solver 
and Matlab's backslash for comparison
are presented in Table \ref{tab2d} 
for different values of $n$. 
The largest numerical experiments in two dimensions 
have been computed for $n = 2049$ with 4 198 401 unknowns.
As can be observed in Table \ref{tab2d}, the CPU times applying Matlab's backslash are about the same 
size as for the initialization, whereas the CPU times of the FFT based fast 
direct
solver grow with 
$\mathcal{O}(N \log N)$.
The observed CPU times are in good agreement with the nearly optimal computational complexity described by Conclusion \ref{conclusion2d} in Section \ref{Sec:FastSolver2D}.

\begin{table}[h!]
\begin{center}
\begin{tabular}{@{}*{1}{c}@{}*{6}{c}@{}} 
  \toprule
   $n$           & 65 & 129 & 257 & 513 & 1025 & 2049 \\
  \midrule
   Initialization                 \hspace{0.5cm} & 0.07  &  0.09  & 0.31 &  1.74 &  13.38  & 118.99 
   \\
   Matlab's backslash     \hspace{0.5cm}  & 0.04 &  0.14  & 0.50 &  2.64  &  12.42 & 93.32 
   \\
   Fast solver                 \hspace{0.5cm}  & 0.01  &  0.02  & 0.06 &  0.23  &  1.05   &  4.58 
   \\
  \bottomrule
\end{tabular}
\end{center}
\caption{CPU times in seconds for different values of $n = n_1 = n_2$ in the two-dimensional 
case}
\label{tab2d}
\end{table}

Table \ref{tab2d_2} presents the CPU times in seconds for various combinations of 
$n_1$ and $n_2$  in the two-dimensional case. More precisely, we compare here the computational times
with respect to a large difference in the magnitude of $n_1$ and $n_2$. Both cases are discussed $n_1 \gg %>> 
n_2$
and $n_2 \gg %>> 
n_1$. As one can observe in Table \ref{tab2d_2} 
it does not influence the computational times of the solver
whether $n_1 \gg %>> 
n_2$ or $n_2 \gg %>> 
n_1$.
For $n_1 = 65$ and $n_2 = 2049$, the fast solver's CPU time was 0.10 seconds, and for
$n_1 = 2049$ and $n_2 = 65$, it was 0.12 seconds.
However, if $n_1$ is very large, it influences the computational times of the initialization process,
more precisely, of
solving the generalized eigenvalue problem 
\eqref{EVP1a} performed
by applying Matlab's function
\texttt{eig}, which requires the full representations of the matrices
$\bsK_1$ and $\bsM_1$, 
which are otherwise stored as sparse matrices.
These increased computational times can be observed in the last two columns of Table
\ref{tab2d_2}. We note here that the solution method described in \cite{HeiRosToi2003}
would be faster for solving this generalized eigenvalue problem.

\begin{table}[h!]
\begin{center}
\begin{tabular}{@{}*{1}{c}@{}*{4}{c}@{}} 
  \toprule
   $n_1$           & 65 & 65     & 2049 & 2049 \\
   $n_2$           & 65 & 2049 & 65 & 2049 \\
  \midrule
   Initialization                 \hspace{0.5cm} & 0.07 & 0.17 & 110.83 & 118.99 \\
   Matlab's backslash     \hspace{0.5cm} & 0.04 & 0.98 & 1.07     &  93.32 \\
   Fast solver                 \hspace{0.5cm}  & 0.01 & 0.10 & 0.12    &  4.58  \\
  \bottomrule
\end{tabular}
\end{center}
\caption{CPU times in seconds for various combinations of $n_1$ and $n_2$ in the two-dimensional 
case}
\label{tab2d_2}
\end{table}

\subsection{Numerical results for the three-dimensional case}
\label{Ssec:Num3d}

In the three-dimensional case, we chose $n_j = n$ for all $j=1,2,3$ for the first set of numerical experiments again.
The CPU times in seconds for the initialization process as well as the application of the fast solver 
and Matlab's backslash for comparison are presented in Table \ref{tab3d} 
for different values of $n$. 
As can be observed from Table \ref{tab3d}, the computational times for applying Matlab's backslash 
rise very quickly. 
%For $n=65$, its CPU time was 608.12 seconds, whereas
%the fast solver's CPU time was only 2.02 seconds.
For $n=65$, its CPU time was 572.91 %608.12 
seconds, whereas
the fast solver's CPU time was only 1.03 %2.02 
seconds.
Then already for $n=129$ with 2 146 689 unknowns,
Matlab runs out of memory when the matrix $\bsA$ is formed which makes it impossible
to solve the problem \eqref{probA} by applying Matlab's backslash.
However, since the FFT based fast solver only needs the setting up 
of the left upper block $\bsC_{bb}$ 
of the matrix $\bsB - \bsA$ and not of the entire matrix $\bsA$,
the solver can be applied solving the problem \eqref{probA}
up to size $n = 513$ with 135 005 697 unknowns (Remark \ref{remarkLaptop}).
Moreover, we can observe from
Table \ref{tab3d} that the CPU times applying Matlab's backslash are already 
much larger for $n=33$
than for the initialization process and the application of the fast solver.
Table \ref{tab3d} shows that the CPU times of applying the FFT based fast 
direct
solver are 
nearly optimal order of complexity $\mathcal{O}(N \log N)$ %.
and match well with the discussed results on the computational complexity in Section \ref{Sec:FastSolver3D}; see Conclusion \ref{conclusion3d}.

\begin{table}[h!]
\begin{center}
\begin{tabular}{@{}*{1}{c}@{}*{7}{c}@{}} 
  \toprule
   $n$           & 9 & 17 & 33 & 65 & 129 & 257 & 513 \\
  \midrule
  Initialization 	 	   \hspace{0.5cm} & 0.07  & 0.07 & 0.08  &   0.45   &  0.41 &  2.16 
  & 17.33  
  \\
   Matlab's backslash    \hspace{0.5cm} & 0.02 & 0.28 & 6.39  & 572.91 &     --     &     --  &  -- \\
   Fast solver 		    \hspace{0.5cm} & 0.09 & 0.12 & 0.21  & 1.03    
   &  8.43 
   & 69.55  
   & 672.27 
   \\
  \bottomrule
\end{tabular}
\end{center}
\caption{CPU times in seconds for different values of $n = n_1 = n_2 = n_3$ in the three-dimensional 
case}
\label{tab3d}
\end{table}

Table \ref{tab3d_2} presents the CPU times in seconds for various combinations of 
$n_j$, $j=1,2,3$, in the three-dimensional 
case similar to Table \ref{tab2d_2} for two dimensions 
comparing computational times
with respect to a large difference in the magnitudes of $n_j$. 
Again we discuss all possible combinations, e.g.,  $n_1 \gg %>> 
n_2, n_3$ or
$n_1, n_2 \gg %>> 
n_3$.
The third and fourth column of Table \ref{tab3d_2}
comparing $n_2 \gg %>> 
n_1, n_3$ and
$n_3 \gg %>> 
n_1, n_2$
show that the CPU times in applying Matlab's backslash are similar (1.06 and 0.96 in 
the third and fourth column, respectively) and also
applying the fast solver 
(0.22 and 0.25 in the third and fourth column, respectively).
However, note that the FFT based direct solver is already 
four times faster than Matlab's backslash for this set of values.
Also, the sixth, seventh and eighth column present similar numerical results,
but Matlab's backslash cannot be applied anymore
for these combinations, since Matlab
runs out of memory when forming the matrix $\bsA$.
Note that the computational times for the initialization process needed
to apply the FFT based direct solver take only around 2 seconds
for these sets of values.

\begin{table}[h!]
\begin{center}
\begin{tabular}{@{}*{1}{c}@{}*{8}{c}@{}} 
  \toprule
   $n_1$           & 9 & 9     & 9 &  513 & 9     & 513  & 513 & 513 \\
   $n_2$           & 9 & 513 & 9 & 9      & 513 & 9      & 513 & 513 \\
   $n_3$           & 9 & 9     & 513 & 9  & 513 &  513 &  9    & 513 \\
  \midrule
   Initialization                 \hspace{0.5cm} & 0.07 & 0.36 & 0.12 & 1.91 &  2.17  & 2.20 & 2.40  & 17.33 \\
   Matlab's backslash     \hspace{0.5cm} & 0.02 &  1.06 & 0.96 & 0.90 & --       & --      & --      & -- \\
   Fast solver                 \hspace{0.5cm}  & 0.09 &  0.22 & 0.25 & 0.62 & 10.95 & 8.00 & 7.01 & 672.27 \\
  \bottomrule
\end{tabular}
\end{center}
\caption{CPU times in seconds for various combinations of $n_1$ and $n_2$ and $n_3$ 
in the three-dimensional case}
\label{tab3d_2}
\end{table}

\begin{remark}
For a large $n$, applying Matlab's LU decomposition \texttt{lu} on the subproblems
\eqref{A2D2}--\eqref{B2D2aux} as well as Matlab's \texttt{kron} function
are the most time consuming functions in the computational process.
Based on this it can be expected that, for example, a \texttt{C++}
implementation would be essentially faster than the Matlab
implementation.
\end{remark}

\section{Conclusions}
\label{Sec:Conclusions}

In this work, we have derived an FFT based direct solver for solving efficiently the Helmholtz equation in a 
rectangular domain with an absorbing boundary condition.
The model problem and fast solver are discussed for two- 
and three-dimensional domains as well as
numerical results for both cases are presented.
%The experiments show nearly optimal order of complexity $\mathcal{O}(N \log N)$
%matching the efficiency of the FFT method.
%
We have shown that the method has the nearly optimal order of complexity $\mathcal{O}(N \log N)$ matching the efficiency of the FFT method. The numerical experiments illustrate that the computational complexity $\mathcal{O}(N \log N)$ is also achievable in practice.
In particular, the numerical results in this work demonstrate %show 
the efficiency of the fast solver
compared to Matlab's backslash and also with respect to the computational memory 
needed to solve the discretized Helmholtz problem for a large number of unknowns.

\section*{Acknowledgments}

The authors gratefully acknowledge the financial support by the
Aca\-demy of Finland
under the grant 295897.
The authors like to thank Dr. Kazufumi Ito for many fruitful
discussions on these fast direct solvers.
Finally, the authors thank the anonymous referee as well as the editor
for the valuable comments improving the paper.

\bibliographystyle{siamplain}
\bibliography{bibliographyToiWol2017arXivNew.bib}

\begin{thebibliography}{10}

\bibitem{BamJolRob1990}
{\sc A.~Bamberger, P.~Joly, and J.~E. Roberts}, {\em Second-order absorbing
  boundary conditions for the wave equation: a solution for the corner
  problem}, SIAM J. Numer. Anal., 27 (1990), pp.~323--352.

\bibitem{Banegas78}
{\sc A.~Banegas}, {\em Fast {P}oisson solvers for problems with sparsity},
  Math. Comp., 32 (1978), pp.~441--446.

\bibitem{BuzDorGeoGol1971}
{\sc B.~Buzbee, F.~W. Dorr, J.~A. George, and G.~H. Golub}, {\em The direct
  solution of the discrete {P}oisson equation on irregular regions}, SIAM J.
  Numer. Anal., 8 (1971), pp.~722--736.

\bibitem{Davis2006}
{\sc T.~A. Davis}, {\em Direct methods for sparse linear systems}, vol.~2,
  Fundamentals of Algorithms, SIAM, Philadelphia, PA, USA, 2006.

\bibitem{ErnGol1994}
{\sc O.~Ernst and G.~H. Golub}, {\em A domain decomposition approach to solving
  the {H}elmholtz equation with a radiation boundary condition}, Contemp.
  Math., 157 (1994), pp.~177--192.

\bibitem{Ern1996}
{\sc O.~G. Ernst}, {\em A finite-element capacitance matrix method for exterior
  {H}elmholtz problems}, Numer. Math., 75 (1996), pp.~175--204.

\bibitem{Guddati04}
{\sc M.~N. Guddati and B.~Yue}, {\em Modified integration rules for reducing
  dispersion in finite element methods}, Comput. Methods Appl. Mech. Engrg.,
  193 (2004), pp.~275--287.

\bibitem{HeiItoToi2018}
{\sc E.~Heikkola, K.~Ito, and J.~Toivanen}, {\em A parallel domain
  decomposition method for the {H}elmholtz equation in layered media},
  submitted for publication,  (2018).

\bibitem{HeiKuzLip1999}
{\sc E.~Heikkola, Y.~A. Kuznetsov, and K.~N. Lipnikov}, {\em Fictitious domain
  methods for the numerical solution of three-dimensional acoustic scattering
  problems}, J. Comput. Acoust., 7 (1999), pp.~161--183.

\bibitem{HeiKuzNeiToi1998}
{\sc E.~Heikkola, Y.~A. Kuznetsov, P.~Neittaanm{\"a}ki, and J.~Toivanen}, {\em
  Fictitious domain methods for the numerical solution of two-dimensional
  scattering problems}, J. Comput. Phys., 145 (1998), pp.~89--109.

\bibitem{HeiRosToi2003}
{\sc E.~Heikkola, T.~Rossi, and J.~Toivanen}, {\em Fast direct solution of the
  {H}elmholtz equation with a perfectly matched layer or an absorbing boundary
  condition}, Internat. J. Numer. Methods Engrg., 57 (2003), pp.~2007--2025.

\bibitem{HeiRosToi2003b}
{\sc E.~Heikkola, T.~Rossi, and J.~Toivanen}, {\em A parallel fictitious domain
  method for the three-dimensional {H}elmholtz equation}, SIAM J. Sci. Comput.,
  24 (2003), pp.~1567--1588.

\bibitem{Ito08}
{\sc K.~Ito, Z.~Qiao, and J.~Toivanen}, {\em A domain decomposition solver for
  acoustic scattering by elastic objects in layered media}, J. Comput. Phys.,
  227 (2008), pp.~8685--8698.

\bibitem{Kuznetsov89}
{\sc Y.~A. Kuznetsov and A.~M. Matsokin}, {\em On partial solution of systems
  of linear algebraic equations}, Sov. J. Numer. Anal. Math. Modelling, 4
  (1989), pp.~453--467.

\bibitem{LynchRiceThomas1964}
{\sc R.~E. Lynch, J.~R. Rice, and D.~H. Thomas}, {\em Direct solution of
  partial difference equations by tensor product methods}, Numerische
  Mathematik, 6 (1964), pp.~185--199.

\bibitem{RosToi1999b}
{\sc T.~Rossi and J.~Toivanen}, {\em A nonstandard cyclic reduction method, its
  variants and stability}, SIAM J. Matrix Anal. Appl., 20 (1999), pp.~628--645.

\bibitem{RosToi1999}
{\sc T.~Rossi and J.~Toivanen}, {\em A parallel fast direct solver for block
  tridiagonal systems with separable matrices of arbitrary dimension}, SIAM J.
  Sci. Comput., 20 (1999), pp.~1778--1793.

\bibitem{Swarztrauber1977}
{\sc P.~N. Swarztrauber}, {\em The methods of cyclic reduction, {F}ourier
  analysis and the {FACR} algorithm for the discrete solution of {P}oisson's
  equation on a rectangle}, SIAM Rev., 19 (1977), pp.~490--501.

\bibitem{TezMacFar2001}
{\sc R.~Tezaur, A.~Macedo, and C.~Farhat}, {\em Iterative solution of
  large-scale acoustic scattering problems with multiple right hand-sides by a
  domain decomposition method with {L}agrange multipliers}, Internat. J. Numer.
  Methods Engrg., 51 (2001), pp.~1175--1193.

\bibitem{Turkel2001}
{\sc E.~Turkel}, {\em Numerical difficulties solving time harmonic systems},
  Nato Sci. S SS III Comput. Systems Sci., 177 (2001), pp.~319--337.

\bibitem{Vas1984}
{\sc P.~Vassilevski}, {\em Fast algorithm for solving a linear algebraic
  problem with separable variables}, Dokl. Bolg. Akad. Nauk., 37 (1984),
  pp.~305--308.

\end{thebibliography}

\end{document}